\documentclass[12 pt, a4 paper]{article}
\usepackage{amssymb,amsmath,graphicx,enumerate}
\usepackage{longtable,mathrsfs}
\usepackage{hyperref}
\numberwithin{equation}{section}
\newtheorem{theorem}{Theorem}[section]
\newtheorem{lemma}[theorem]{Lemma}
\newtheorem{definition}{Definition}[section]
\newtheorem{remark}{Remark}[section]

\allowdisplaybreaks

\begin{document}
\begin{center}\begin{large}
Some Minkowski-type inequalities using generalized proportional Hadamard  fractional integral operators
\end{large}\end{center}
\begin{center}
                 $ Asha \, B. Nale^{1},\, Satish  \, K. Panchal ^{2}, \, Vaijanath  \, L. Chinchane ^{3}$

$^{1,2}$Department of Mathematics,\\
Dr. Babasaheb Ambedkar Marathwada\\
University, Aurangabad-431 004, INDIA.\\
$^{3}$Department of Mathematics,\\
Deogiri Institute of Engineering and Management\\
Studies Aurangabad-431005, INDIA\\
ashabnale@gmail.com/drskpanchal@gmail.com/chinchane85@gmail.com
\end{center}
\begin{abstract}
The main objective of  present investigation to obtain some Minkowski-type fractional integral inequalities using generalised proportional Hadamard  fractional integral operators which is introduced by Rahman et al in the paper (certain inequalities via generalized proportional Hadamard fractional integral operators), Advances in Differential Equations, (2019), 2019:454. In addition, we establish some other fractional integral inequalities for positive and continuous functions.
\end{abstract}
\section{Introduction}
 \paragraph{}Fractional calculus is generalization of traditional calculus into non-integer differential and integral order. Fractional calculus is very important due to it's various application in field of science and technology. Fractional integral inequalities play big role in obtaining uniqueness of solution of fractional ordinary differential equations, fractional partial differential equations and fractional Boundary value problems. Recently, a number of researchers in the fields of fractional integral inequalities  have established different integral inequalities about fractional integral operators such as Riemann-Liouville, Hadamard, Saigo, generalized Katugamapola, Erd$\acute{e}$lyi-Kober, Riemann-Liouville k-fractional, Hadamard k-fractional, (k,s)-Riemann-Liouville and  k-generalized (in terms hypergeometric function) fractional integral operators, see \cite{A1,C2,C3,C4,C5,D2,D3,D4,D5,D6,K1,K2,N1,P1,P2,P3,S1,S4,V1}. In \cite{D1}, Dahmani investigated reverse Minkowski fractional integral inequality by employing Riemann-Liouville  fractional integral . Ahmed Anber and et al \cite{A2} presented some  fractional integral inequalities which is simillar to Minkowski fractional integral inequality  using Riemann-Liouville  fractional integral.  V. L. Chinchane et al \cite{C1} proposed fractional inequalities similar to Minkowski type via Saigo fractional integral operator. In \cite{C6}, V. L. Chinchane studied the reverse Minkowski fractional integral inequality by considering generalized K-fractional integral operator is in terms of the Gauss herpergeometric function. S. Mubeen and et al \cite{M1} have introduced Minkowski inequality involving generalized k-fractional conformable integrals.\\
 G. Rahman et al. \cite{R1,R2,R3} established Minkowski inequality and some other fractional inequalities for convex functions by employing fractional proportional integral operators. In \cite{J1,R4,R5}, F. Jarad et al and  G.Rahman presented concepts of non-local fractional proportional and generalized Hadamard proportional integrals involving exponential functions in their kernels.
\paragraph{} Motivated from \cite{R2,R3,R4,R5}, our purpose in this paper is to propose some new results using generalized Hadamard proportional integrals. The paper has been organized as follows, in Section 2, we recall basic definitions, remarks and lemma related to generalized Hadamard proportional integrals. In Section 3, we obtain reverse Minkowski fractional integral inequality using generalized Hadamard proportional integrals, In Section 4, we present some other inequalities using generalized Hadamard proportional integrals.

\section{preliminary}
Here, we present some important definition, remarks and lemma of generalised proportional Hadamard fractional integral operator which will be used throughout this paper.
\begin{definition} The left and right sided generalized proportional fractional integrals are respectively defined by

\begin{equation}
(_a\mathfrak{J}^{\alpha,\beta}z)(x)=\frac{1}{\beta^{\alpha}\Gamma(\alpha)} \int_{a}^{x}e^{[\frac{\beta-1}{\beta}(x-t)]}(x-t)^{\alpha-1}z(t)dt, a<x,
\end{equation}
and
\begin{equation}
(\mathfrak{J}^{\alpha,\beta}_{b}z)(x)=\frac{1}{\beta^{\alpha}\Gamma(\alpha)} \int_{x}^{b}e^{[\frac{\beta-1}{\beta}(t-x)]}(t-x)^{\alpha-1}z(t)dt, x<b,
\end{equation}
where the proportionality index $\beta \in (0,1]$ and $\alpha \in \mathbb{C}$ with $\mathfrak{R(\alpha)}> 0$.
\end{definition}
\begin{remark}
If we consider $\beta=1$ in (2.1) and (2.2), then we get the well known left and right Riemann-Liouville integrals which are respectively defined by
\begin{equation}
(_a\mathfrak{J}^{\alpha}z)(x)=\frac{1}{\Gamma(\alpha)} \int_{a}^{x}(x-t)^{\alpha-1}z(t)dt, a<x,
\end{equation}
and
\begin{equation}
(\mathfrak{J}^{\alpha}_{b}z)(x)=\frac{1}{\Gamma(\alpha)} \int_{x}^{b}(t-x)^{\alpha-1}z(t)dt, x<b,
\end{equation}
where $\alpha \in \mathbb{C}$ with $\mathfrak{R(\alpha)}>0 $.
\end{remark}
Recently, Rahman et al.\cite{R5} proposed the following generalized Hadamard proportional fractional integrals.

\begin{definition} The left sided generalized Hadamard proportional fractional integral of order $\alpha > 0$ and proportional index $\beta \in(0,1]$ is
defined by
\begin{equation}
(_a\mathcal{H}^{\alpha,\beta}z)(x)=\frac{1}{\beta^{\alpha}\Gamma(\alpha)} \int_{a}^{x}e^{[\frac{\beta-1}{\beta}(ln x-ln t)]}(ln x-ln t)^{\alpha-1}\frac{z(t)}{t}dt, a<x.
\end{equation}
\end{definition}

\begin{definition} The right sided generalized Hadamard proportional fractional integral of order $\alpha > 0$ and proportional index $\beta \in(0,1]$ is defined by
\begin{equation}
(\mathcal{H}^{\alpha,\beta}_{b}z)(x)=\frac{1}{\beta^{\alpha}\Gamma(\alpha)} \int_{x}^{b}e^{[\frac{\beta-1}{\beta}(ln t-ln x)]}(ln t-ln x)^{\alpha-1}\frac{z(t)}{t}dt, x<b.
\end{equation}
\end{definition}
\begin{definition} The one sided generalized Hadamard proportional fractional integral of order $\alpha > 0$ and proportional index $\beta \in(0,1]$ is
defined by
\begin{equation}
(\mathcal{H}^{\alpha,\beta}_{1,x}z)(x)=\frac{1}{\beta^{\alpha}\Gamma(\alpha)} \int_{1}^{x}e^{[\frac{\beta-1}{\beta}(ln x-ln t)]}(ln x-ln t)^{\alpha-1}\frac{z(t)}{t}dt, t>1,
\end{equation}
where $\Gamma(\alpha)$ is the classical well known gamma function.
\end{definition}
\begin{remark} If we consider $\beta=1$, then (2.5)-(2.7) will led to the following well known Hadamard fractional integrals
\begin{equation}
(_a\mathcal{H}^{\alpha}z)(x)=\frac{1}{\Gamma(\alpha)} \int_{a}^{x}(ln x-ln t)^{\alpha-1}\frac{z(t)}{t}dt, a<x,
\end{equation}

\begin{equation}
(\mathcal{H}^{\alpha}_{b}z)(x)=\frac{1}{\Gamma(\alpha)} \int_{x}^{b}(ln t-ln x)^{\alpha-1}\frac{z(t)}{t}dt, x<b,
\end{equation}

and
\begin{equation}
(\mathcal{H}^{\alpha}_{1,x}z)(x)=\frac{1}{\Gamma(\alpha)} \int_{1}^{x}(ln x-ln t)^{\alpha-1}\frac{z(t)}{t}dt, x>1.
\end{equation}
\end{remark}
One can easily prove the following results:
\begin{lemma}

\begin{equation}
(\mathcal{H}^{\alpha,\beta}_{1,x}e^{[\frac{\beta-1}{\beta}(ln x)]}(ln x)^{\lambda-1})(x)=\frac{\Gamma(\lambda)}{\beta^{\alpha}\Gamma(\alpha+\lambda)}e^{[\frac{\beta-1}{\beta}(ln x)]}(ln x)^{\alpha+\lambda-1},
\end{equation}
and the semigroup property
\begin{equation}
(\mathcal{H}^{\alpha,\beta}_{1,x})(\mathcal{H}^{\lambda,\beta}_{1,x})z(x)=(\mathcal{H}^{\alpha+\lambda,\beta}_{1,x})z(x).
\end{equation}
\end{lemma}
\begin{remark} If $\beta =1$, then (2.11) will reduce to the result of \cite{S2} as defined by
\begin{equation}
(\mathcal{H}^{\alpha}_{1,x}(ln x)^{\lambda-1})(x)=\frac{\Gamma(\lambda)}{\Gamma(\alpha+\lambda)}(ln x)^{\alpha+\lambda-1}.
\end{equation}
\end{remark}
 \section{Reverse Minkowski fractional integral inequality }
In this section, we establish reverse Minkowski fractional integral inequality involving generalized proportional Hadamard  fractional integral operators.
 \begin{theorem} Let $p\geq1$ and let $f$, $g$ be two positive function on $[1, \infty)$, such that for all $x>1$, $\mathcal{H}^{\alpha,\beta}_{1,x}[f^{p}(x)]<\infty$, $\mathcal{H}^{\alpha,\beta}_{1,x}[g^{q}(x)]<\infty$. If $0<m\leq \frac{f(\tau)}{g(\tau)}\leq M$, $\tau \in (1,x)$ we have
 \begin{equation}
\left[\mathcal{H}^{\alpha,\beta}_{1,x}[f^{p}(x)]\right]^{\frac{1}{p}}+\left[\mathcal{H}^{\alpha,\beta}_{1,x}[g^{q}(x)]\right]^{\frac{1}{p}}\leq \frac{1+M(m+2)}{(m+1)(M+1)}\left[\mathcal{H}^{\alpha,\beta}_{1,x}[(f+g)^{p}(x)]\right]^{\frac{1}{p}},
\end{equation}
where $\alpha>0,$ $\beta \in (0,1]$ $\alpha \in \mathbb{C}$ and $\mathfrak{R(\alpha)}> 0.$
\end{theorem}
\textbf{Proof}: Using the condition $\frac{f(\tau)}{g(\tau)}\leq M$, $\tau \in (1,x)$, $x>1$, we can write
\begin{equation}
(M+1)^{p}f(\tau)\leq M^{p}(f+g)^{p}(\tau).
\end{equation}
Consider\\
\begin{equation}
\mathrm{\psi}(x,\tau)=\frac{1}{\beta^{\alpha}\Gamma(\alpha)\tau}e^{[\frac{\beta-1}{\beta}(ln x-ln \tau)]}(ln x-ln \tau)^{\alpha-1}.
\end{equation}
Clearly, we can say that the function $\mathrm{\psi}(x,\tau)$ remain positive because for all $\tau \in (1, x)$ , $(x>1), \alpha, \beta>0$.
Multiplying both side of (3.2) by $\mathrm{\psi}(x,\tau)$, then integrating resulting identity with respect to $\tau$ from $1$ to $x$, we get
\begin{equation}
\begin{split}
&\frac{(M+1)^{p}}{\beta^{\alpha}\Gamma(\alpha)} \int_{1}^{x}e^{[\frac{\beta-1}{\beta}(ln x-ln \tau)]}(ln x-ln \tau)^{\alpha-1}f^{p}(\tau)\frac{d\tau}{\tau}\\
&\leq \frac{M^{p}}{\beta^{\alpha}\Gamma(\alpha)} \int_{1}^{x}e^{[\frac{\beta-1}{\beta}(ln x-ln \tau)]}(ln x-ln \tau)^{\alpha-1}(f+g)^{p}(\tau)\frac{d\tau}{\tau},
\end{split}
\end{equation}
\noindent which is equivalent to
\begin{equation}
\mathcal{H}^{\alpha,\beta}_{1,x}[f^{p}(x)]  \leq \frac{M^{p}}{(M+1)^{p}} \left[\mathcal{H}^{\alpha,\beta}_{1,x}[(f+g)^{p}(x)]\right],
\end{equation}
\noindent hence, we can write
\begin{equation}
\left[\mathcal{H}^{\alpha,\beta}_{1,x}[f^{p}(x)] \right]^{\frac{1}{p}} \leq \frac{M}{(M+1)} \left[\mathcal{H}^{\alpha,\beta}_{1,x}[(f+g)^{p}(x)]\right]^{\frac{1}{p}}.
\end{equation}
On other hand, using condition $m\leq \frac{f(\tau)}{g(\tau)}$, we obtain
\begin{equation}
(1+\frac{1}{m})g(\tau)\leq \frac{1}{m}(f(\tau)+g(\tau)),
\end{equation}
therefore,
\begin{equation}
(1+\frac{1}{m})^{p}g^{p}(\tau)\leq(\frac{1}{m})^{p}(f(\tau)+g(\tau))^{p}.
\end{equation}
Now, multiplying both side of (3.8) by  $\mathrm{\psi}(x,\tau)$, ( $\tau \in(1,x)$, $x>1$), where  $G(x,\tau)$ is defined by (3.3). Then integrating resulting identity with respect to $\tau$ from $1$ to $x$, we have
\begin{equation}
\left[\mathcal{H}^{\alpha,\beta}_{1,x}[g^{p}(x)]\right]^{\frac{1}{p}}  \leq \frac{1}{(m+1)} \left[\mathcal{H}^{\alpha,\beta}_{1,x}[(f+g)^{p}(x)]\right]^{\frac{1}{p}}.
\end{equation}
The inequalities (3.1) follows on adding the inequalities (3.6) and (3.9).\\
Our second result is as follows.
\begin{theorem} Let $p\geq1$ and $f$, $g$ be two positive function on $[1, \infty)$, such that for all $x>0$, $\mathcal{H}^{\alpha,\beta}_{1,x}[f^{p}(x)]<\infty$, $\mathcal{H}^{\alpha,\beta}_{1,x}[g^{q}(x)]<\infty$. If $0<m\leq \frac{f(\tau)}{g(\tau)}\leq M$, we have
\begin{equation}
\begin{split}
\left[\mathcal{H}^{\alpha,\beta}_{1,x}[f^{p}(x)] \right]^{\frac{2}{p}}+\left[\mathcal{H}^{\alpha,\beta}_{1,x}[g^{q}(x)] \right]^{\frac{2}{p}}\geq &(\frac{(M+1)(m+1)}{M}-2)\left[\mathcal{H}^{\alpha,\beta}_{1,x}[f^{p}(x)] \right]^{\frac{1}{p}}+\\
&\left[I^{\alpha, \beta, \eta}_{0, x}[g^{q}(x)] \right]^{\frac{1}{p}},
\end{split}
\end{equation}
where $\alpha>0,$ $\beta \in (0,1]$ $\alpha \in \mathbb{C}$ and $\mathfrak{R(\alpha)}> 0.$
\end{theorem}
\textbf{Proof}: Multiplying the inequalities (3.6) and (3.9), we obtain
\begin{equation}
\frac{(M+1)(m+1)}{M}\left[\mathcal{H}^{\alpha,\beta}_{1,x}[f^{p}(x)]\right]^{\frac{1}{p}}\, \left[\mathcal{H}^{\alpha,\beta}_{1,x}[g^{q}(x)]\right]^{\frac{1}{p}}\leq \left([\mathcal{H}^{\alpha,\beta}_{1,x}[(f(x)+g(x))^{p}]]^{\frac{1}{p}}\right)^{2}.
\end{equation}
Applying Minkowski inequalities to the right hand side of (3.11), we have
 \begin{equation}
 (\left[\mathcal{H}^{\alpha,\beta}_{1,x}[(f(x)+g(x))^{p}]\right]^{\frac{1}{p}})^{2}\leq (\left[\mathcal{H}^{\alpha,\beta}_{1,x}[f^{p}(x)]\right]^{\frac{1}{p}}+\left[\mathcal{H}^{\alpha,\beta}_{1,x}[g^{q}(x)]\right]^{\frac{1}{p}})^{2},
\end{equation}
which implies that
\begin{equation}
\begin{split}
 \left[\mathcal{H}^{\alpha,\beta}_{1,x}[(f(x)+g(x))^{p}]\right]^{\frac{2}{p}}\leq & \left[\mathcal{H}^{\alpha,\beta}_{1,x}[f^{p}(x)]\right]^{\frac{2}{p}}+
 \left[\mathcal{H}^{\alpha,\beta}_{1,x}[g^{q}(x)]\right]^{\frac{2}{p}}\\
 &+2\left[\mathcal{H}^{\alpha,\beta}_{1,x}[f^{p}(x)]\right]^{\frac{1}{p}} \left[\mathcal{H}^{\alpha,\beta}_{1,x}[g^{q}(x)]\right]^{\frac{1}{p}},
\end{split}
\end{equation}
using (3.11) and (3.13) we obtain (3.10).
Theorem 3.2 is thus proved.
\section{ Fractional integral inequalities involving Saigo fractional integral operator}
Here, we establish some new integral inequalities using  generalized proportional Hadamard  fractional integral operators.
\begin{theorem} Let $p>1$,  $\frac{1}{p}+\frac{1}{q}=1 $ and $f$, $g$ be two positive function on $[1, \infty)$, such that $\mathcal{H}^{\alpha,\beta}_{1,x}[f(x)]<\infty$, $\mathcal{H}^{\alpha,\beta}_{1,x}[g(x)]<\infty$. If $0<m\leq \frac{f(\tau)}{g(\tau)}\leq M < \infty$, $\tau \in (1,x)$ we have
\begin{equation}
\left[\mathcal{H}^{\alpha,\beta}_{1,x}[f(x)]\right]^{\frac{1}{p}} \left[\mathcal{H}^{\alpha,\beta}_{1,x}[g(x)]\right]^{\frac{1}{q}}
\leq (\frac{M}{m})^{\frac{1}{pq}}\left[\mathcal{H}^{\alpha,\beta}_{1,x}[[f(x)]^{\frac{1}{p}}[g(t)]^{\frac{1}{q}}]\right],
\end{equation}
hold, where $\alpha>0,$ $\beta \in (0,1]$ $\alpha \in \mathbb{C}$ and $\mathfrak{R(\alpha)}> 0.$

\end{theorem}
\textbf{Proof:-} Since $\frac{f(\tau)}{g(\tau)}\leq M $, $\tau \in[1,x]$  $x> 1$, therefore \\
\begin{equation}
[g(\tau)]^{\frac{1}{p}}\geq M^{\frac{-1}{q}}[f(\tau)]^{\frac{1}{q}},
\end{equation}
and also,
\begin{equation}
\begin{split}
[f(\tau)]^{\frac{1}{p}}[g(\tau)]^{\frac{1}{q}}&\geq M^{\frac{-1}{q}}[f(\tau)]^{\frac{1}{q}}[f(\tau)]^{\frac{1}{p}}\\
&\geq  M^{\frac{-1}{q}}[f(\tau)]^{\frac{1}{q}+\frac{1}{q}}\\
&\geq  M^{\frac{-1}{q}}[f(\tau)].
\end{split}
\end{equation}
Multiplying both side of (4.3) by $\mathrm{\psi}(x,\tau)$, ( $\tau \in(1,x)$, $x>1$), where  $\mathrm{\psi}(x,\tau)$ is defined by (3.3). Then integrating resulting identity with respect to $\tau$ from $1$ to $x$, we have
\begin{equation}
\begin{split}
&\frac{1}{\beta^{\alpha}\Gamma(\alpha)} \int_{1}^{x}e^{[\frac{\beta-1}{\beta}(ln x-ln \tau)]}(ln x-ln \tau)^{\alpha-1}f(\tau)^{\frac{1}{p}}g(\tau)^{\frac{1}{q}}(\tau)\frac{d\tau}{\tau} \\
&\leq \frac{ M^{\frac{-1}{q}}}{\beta^{\alpha}\Gamma(\alpha)} \int_{1}^{x}e^{[\frac{\beta-1}{\beta}(ln x-ln \tau)]}(ln x-ln \tau)^{\alpha-1}f(\tau))\frac{d\tau}{\tau},
\end{split}
\end{equation}
which implies that,
\begin{equation}
\mathcal{H}^{\alpha,\beta}_{1,x}\left[[f(x)]^{\frac{1}{p}}[g(x)]^{\frac{1}{q}} \right] \leq M^{\frac{-1}{q}} \left[\mathcal{H}^{\alpha,\beta}_{1,x}f(x)\right].
\end{equation}
Consequently,\\
\begin{equation}
\left(\mathcal{H}^{\alpha,\beta}_{1,x}\left[[f(x)]^{\frac{1}{p}}[g(x)]^{\frac{1}{q}} \right]\right)^{\frac{1}{p}} \leq M^{\frac{-1}{pq}} \left[\mathcal{H}^{\alpha,\beta}_{1,x}f(x)\right]^{\frac{1}{p}},
\end{equation}
on other hand, since $m g(\tau)\leq f(\tau)$, \, $\tau \in[1,x)$, $x>1$, then we have
\begin{equation}
[f(\tau)]^{\frac{1}{p}}\geq m^{\frac{1}{p}}[g(\tau)]^{\frac{1}{p}},
\end{equation}
multiplying equation (4.7) by $[g(\tau)]^{\frac{1}{q}}$, we have,\\
\begin{equation}
[f(\tau)]^{\frac{1}{p}}[g(\tau)]^{\frac{1}{q}}\geq m^{\frac{1}{p}}[g(\tau)]^{\frac{1}{q}}[g(\tau)]^{\frac{1}{p}}= m^{\frac{1}{p}}[g(\tau)].
\end{equation}
Multiplying both side of (4.8) by $\mathrm{\psi}(x,\tau)$, ( $\tau \in(1,x)$, $x>1$), where $\mathrm{\psi}(x,\tau)$ is defined by (3.3). Then integrating resulting identity with respect to $\tau$ from $1$ to $x$, we have
\begin{equation}
\begin{split}
&\frac{1}{\beta^{\alpha}\Gamma(\alpha)} \int_{1}^{x}e^{[\frac{\beta-1}{\beta}(ln x-ln \tau)]}(ln x-ln \tau)^{\alpha-1}f(\tau)^{\frac{1}{p}}g(\tau)^{\frac{1}{q}}(\tau)\frac{d\tau}{\tau} \\
&\leq \frac{m^{\frac{1}{p}}}{\beta^{\alpha}\Gamma(\alpha)} \int_{1}^{x}e^{[\frac{\beta-1}{\beta}(ln x-ln \tau)]}(ln x-ln \tau)^{\alpha-1}[g(\tau)]\frac{d\tau}{\tau},
\end{split}
\end{equation}
that is,
\begin{equation}
\mathcal{H}^{\alpha,\beta}_{1,x}\left[[f(x)]^{\frac{1}{p}}[g(x)]^{\frac{1}{q}} \right] \leq m^{\frac{1}{p}} \left[\mathcal{H}^{\alpha,\beta}_{1,x}g(x)\right].
\end{equation}
Hence we can write,
\begin{equation}
\left(\mathcal{H}^{\alpha,\beta}_{1,x}\left[[f(x)]^{\frac{1}{p}}[g(x)]^{\frac{1}{q}} \right]\right)^{\frac{1}{q}} \leq m^{\frac{1}{pq}} \left[\mathcal{H}^{\alpha,\beta}_{1,x}f(x)\right]^{\frac{1}{q}},
\end{equation}
multiplying equation (4.6) and (4.11) we get the result (4.1).
\begin{theorem} Let $f$ and $g$ be two positive function on $[1, \infty[$, such that\\  $\mathcal{H}^{\alpha,\beta}_{1,x}[f^{p}(x)]<\infty$,
 $\mathcal{H}^{\alpha,\beta}_{1,x}[g^{q}(x)]<\infty$. $x>1$,  If $0<m\leq \frac{f(\tau)^{p}}{g(\tau)^{q}}\leq M < \infty$, $\tau \in [1,x]$. Then we have
 \begin{equation}
\left[\mathcal{H}^{\alpha,\beta}_{1,x}f^{p}(x)\right]^{\frac{1}{p}} \left[\mathcal{H}^{\alpha,\beta}_{1,x}g^{q}(x)\right]^{\frac{1}{q}}\leq (\frac{M}{m})^{\frac{1}{pq}}\left[\mathcal{H}^{\alpha,\beta}_{1,x}(f(x)g(x))\right] hold,
\end{equation}
where $p>1$, $\frac{1}{p}+\frac{1}{q}=1 $,  $\alpha>0,$ $\beta \in (0,1]$ $\alpha \in \mathbb{C}$ and $\mathfrak{R(\alpha)}> 0.$

\end{theorem}
\textbf{Proof:-}
Replacing $f(\tau)$ and $g(\tau)$ by $f(\tau)^{p}$ and $g(\tau)^{q}$,  $\tau \in [1,x]$, $x>1$ in theorem 4.1, we obtain (4.12).
\paragraph{} Now, here we obtain some fractional integral inequality related to Minkowsky inequality as follows
\begin{theorem} let $f$ and $g$ be two integrable functions on $[1, \infty]$ such that $\frac{1}{p}+\frac{1}{q}=1, p>1,$ and $0<m<\frac{f(\tau)}{g(\tau)}<M, \tau \in (1,x), x>1.$ Then we have
\begin{equation}
\mathcal{H}^{\alpha,\beta}_{1,x}\{fg\}(x)\leq \frac{2^{p-1}M^{p}}{p(M+1)^{p}}\left(\mathcal{H}^{\alpha,\beta}_{1,x}[f^{p}+g^{p}](x)\right)+\frac{2^{q-1}}{q(m+1)^{q}}\left(\mathcal{H}^{\alpha,\beta}_{1,x}[f^{q}+g^{q}](x)\right),
\end{equation}
where $\alpha>0,$ $\beta \in (0,1]$ $\alpha \in \mathbb{C}$ and $\mathfrak{R(\alpha)}> 0.$
\end{theorem}
\textbf{Proof:-} Since, $\frac{f(\tau)}{g(\tau)}<M, \tau \in (1,x), x>1,$ we have
\begin{equation}
(M+1)f(\tau)\leq M(f+g)(\tau).
\end{equation}
Taking $p^{th}$  power on both side of (4.14)and multiplying resulting identity by $ \mathrm{\psi}(x,\tau)$, which is positive because $\tau \in [1,x]$, $x>1$, then integrate the resulting identity with respect to $\tau$ from $1$ to $x$, we get
\begin{equation}
\begin{split}
&\frac{(M+1)^{p}}{\beta^{\alpha}\Gamma(\alpha)} \int_{1}^{x}e^{[\frac{\beta-1}{\beta}(ln x-ln \tau)]}(ln x-ln \tau)^{\alpha-1}f(\tau)^{p}\frac{d\tau}{\tau}\\
 &\leq \frac{M^{p}}{\beta^{\alpha}\Gamma(\alpha)} \int_{1}^{x}e^{[\frac{\beta-1}{\beta}(ln x-ln \tau)]}(ln x-ln \tau)^{\alpha-1}[(f+g)^{p}(\tau)]\frac{d\tau}{\tau},
\end{split}
\end{equation}
therefore,
\begin{equation}
\mathcal{H}^{\alpha,\beta}_{1,x}[f^{p}(x)]\leq \frac{M^{p}}{(M+1)^{p}}\mathcal{H}^{\alpha,\beta}_{1,x}[(f+g)^{p}(x)],
\end{equation}
on other hand, $0<m<\frac{f(\tau)}{g(\tau)}, \tau \in (1,x), x>1,$ we can write
\begin{equation}
(m+1)g(\tau)\leq (f+g)(\tau),
\end{equation}
Taking $q^{th}$  power on both side (4.17) and multiplying resulting identity by $ \mathrm{\psi}(x,\tau)$, which is positive because $\tau \in [1,x]$, $x>1$ then integrate the resulting identity with respect to $\tau$ from $1$ to $x$, we get

\begin{equation}
\begin{split}
&\frac{(m+1)^{q}}{\beta^{\alpha}\Gamma(\alpha)} \int_{1}^{x}e^{[\frac{\beta-1}{\beta}(ln x-ln \tau)]}(ln x-ln \tau)^{\alpha-1}g^{q}(\tau)\frac{d\tau}{\tau}\\
 &\leq \frac{1}{\beta^{\alpha}\Gamma(\alpha)} \int_{1}^{x}e^{[\frac{\beta-1}{\beta}(ln x-ln \tau)]}(ln x-ln \tau)^{\alpha-1}[(f+g)^{q}(\tau)]\frac{d\tau}{\tau},
 \end{split}
\end{equation}
consequently, we have
\begin{equation}
\mathcal{H}^{\alpha,\beta}_{1,x}[g^{q}(x)]\leq \frac{1}{(m+1)^{q}}\mathcal{H}^{\alpha,\beta}_{1,x}[(f+g)^{q}(x)].
\end{equation}
Now, using Young inequality
\begin{equation}
[f(\tau)g(\tau)]\leq \frac{f^{p}(\tau)}{p}+\frac{g^{q}(\tau)}{q}.
\end{equation}
Multiplying both side of (4.20) by $ \mathrm{\psi}(x,\tau)$, which is positive because $\tau \in [1,x]$, $x>1$, then integrate the resulting identity with respect to $\tau$ from $1$ to $x$, we get
\begin{equation}
\mathcal{H}^{\alpha,\beta}_{1,x}[f(x)g(x))]\leq \frac{1}{p}\,\mathcal{H}^{\alpha,\beta}_{1,x}[f^{p}(x)]+\frac{1}{q}\,\mathcal{H}^{\alpha,\beta}_{1,x}[g^{q}(x)],
\end{equation}
from equation (4.16), (4.19) and (4.21) we get
\begin{equation}
\mathcal{H}^{\alpha,\beta}_{1,x}[f(x)g(x))]\leq \frac{M^{p}}{p(M+1)^{p}}\,\mathcal{H}^{\alpha,\beta}_{1,x}[(f+g)^{p}(x)]+\frac{1}{q(m+1)^{q}}\,\mathcal{H}^{\alpha,\beta}_{1,x}[(f+g)^{q}(x)],
\end{equation}
now using the inequality $(a+b)^{r}\leq 2^{r-1}(a^{r}+b^{r}), r>1, a,b \geq 0,$ we have
\begin{equation}
\mathcal{H}^{\alpha,\beta}_{1,x}[(f+g)^{p}(x)] \leq 2^{p-1}\mathcal{H}^{\alpha,\beta}_{1,x}[(f^{p}+g^{p})(x)],
\end{equation}
and
\begin{equation}
\mathcal{H}^{\alpha,\beta}_{1,x}[(f+g)^{q}(x)] \leq 2^{q-1}\mathcal{H}^{\alpha,\beta}_{1,x}[(f^{q}+g^{q})(x)].
\end{equation}
Injecting (4.23), (2.24) in (4.22) we get required inequality (4.13). This complete the proof.

\begin{theorem} Let $f$, $g$ be two positive functions defined on $[1, \infty)$, such that  $g$ is non-decreasing. If
\begin{equation}
\mathcal{H}^{\alpha,\beta}_{1,x}f(x)\geq \mathcal{H}^{\alpha,\beta}_{1,x}g(x), x>1.
\end{equation}
then for all $\alpha>0,$ $\beta \in (0,1]$ $\alpha \in \mathbb{C}$ and $\mathfrak{R(\alpha)}> 0$,
\begin{equation}
\mathcal{H}^{\alpha,\beta}_{1,x}f^{\gamma-\delta}(x)\leq  \mathcal{H}^{\alpha,\beta}_{1,x}f^{\gamma}(x)g^{-\delta}(x),
\end{equation}
hold.
\end{theorem}
\textbf{Proof:-}We use arithmetic-geometric inequality, for $\gamma>0$, $\delta>0$, we have:
 \begin{equation}
\frac{\gamma}{\gamma-\delta}f^{\gamma-\delta}(\tau)-\frac{\delta}{\gamma-\delta}g^{\gamma-\delta}(\tau)\leq f^{\gamma}(\tau)g^{-\delta}(\tau), \, \, \tau \in(1,x), x>1.
\end{equation}
Now, multiplying both side of (4.27)  by  $\mathrm{\psi}(x,\tau)$, which is positive because $\tau \in [1,x]$, $x>1$, then integrate the resulting identity with respect to $\tau$ from $1$ to $x$, we get

\begin{equation}
\begin{split}
&\frac{\gamma}{\gamma-\delta\beta^{\alpha}\Gamma(\alpha)} \int_{1}^{x}e^{[\frac{\beta-1}{\beta}(ln x-ln \tau)]}(ln x-ln \tau)^{\alpha-1}f^{\gamma-\delta}(\tau)\frac{d\tau}{\tau}\\
&- \frac{\delta}{\gamma-\delta \beta^{\alpha}\Gamma(\alpha)} \int_{1}^{x}e^{[\frac{\beta-1}{\beta}(ln x-ln \tau)]}(ln x-ln \tau)^{\alpha-1}g^{\gamma-\delta}(\tau)\frac{d\tau}{\tau} \\
&\leq \frac{1}{\beta^{\alpha}\Gamma(\alpha)} \int_{1}^{x}e^{[\frac{\beta-1}{\beta}(ln x-ln \tau)]}(ln x-ln \tau)^{\alpha-1}f^{\gamma}(\tau)g^{-\delta}(\tau)\frac{d\tau}{\tau},
\end{split}
\end{equation}
consequently,
\begin{equation}
\frac{\gamma}{\gamma-\delta}\,\mathcal{H}^{\alpha,\beta}_{1,x}[f^{\gamma-\delta}(x)]-\frac{\delta}{\gamma-\delta} \mathcal{H}^{\alpha,\beta}_{1,x}[g^{\gamma-\delta}(x)] \leq  \mathcal{H}^{\alpha,\beta}_{1,x}[f^{\gamma}(x)g^{-\delta}(x)],
\end{equation}
which implies that,
\begin{equation}
\frac{\gamma}{\gamma-\delta}\,\mathcal{H}^{\alpha,\beta}_{1,x}[f^{\gamma-\delta}(x)]
 \leq \,  \mathcal{H}^{\alpha,\beta}_{1,x}[f^{\gamma}(x)g^{-\delta}(x)]+\frac{\delta}{\gamma-\delta}  \, \mathcal{H}^{\alpha,\beta}_{1,x}[g^{\gamma-\delta}(x)],
\end{equation}
that is
\begin{equation}
 \mathcal{H}^{\alpha,\beta}_{1,x}[f^{\gamma-\delta}(x)]
 \leq \, \frac{\gamma-\delta}{\gamma}\, \mathcal{H}^{\alpha,\beta}_{1,x}[f^{\gamma}(t)g^{-\delta}(x)]+\frac{\delta}{\gamma} \,\mathcal{H}^{\alpha,\beta}_{1,x}[f^{\gamma-\delta}(x)],
\end{equation}
thus we get the result (4.26).
 \begin{theorem} Suppose that $f$, $g$ and $h$ be positive and continuous functions on $[1,\infty)$, such that
 \begin{equation}
(g(\tau)-g(\sigma)) \left(\frac{f(\sigma)}{h(\sigma)}-\frac{f(\tau)}{h(\tau)}\right)\geq 0; \  \tau, \sigma \in[1,x)\,  x>1,
\end{equation}
then for all $\alpha>0,$ $\beta \in (0,1]$ $\alpha \in \mathbb{C}$ and $\mathfrak{R(\alpha)}> 0$,
\begin{equation}
\frac{\mathcal{H}^{\alpha,\beta}_{1,x}[f(x)]}{\mathcal{H}^{\alpha,\beta}_{1,x}[h(x)]}\geq \frac{\mathcal{H}^{\alpha,\beta}_{1,x}[(g f)(x)]}{\mathcal{H}^{\alpha,\beta}_{1,x}[(gh)(x)]}.
\end{equation}
Hold.
\end{theorem}
\textbf{Proof}: Since $f$, $g$ and $h$ be three  positive and continuous functions on $[1,\infty[$ by (4.32), we can write
  \begin{equation}
 g(\tau)\frac{f(\sigma)}{h(\sigma)}+g(\sigma)\frac{f(\tau)}{h(\tau)}-g(\sigma)\frac{f(\sigma)}{h(\sigma)}-g(\tau)\frac{f(\tau)}{h(\tau)}\geq 0; \ \ \tau, \sigma \in(1,x),\ \  x>1.
\end{equation}
 \noindent  Now, multiplying equation (4.34) by $ h(\sigma) h(\tau)$, on both side, we have,
\begin{equation}
g(\tau)f(\sigma)h(\tau)-g(\tau)f(\tau)h(\sigma)-g(\sigma)f(\sigma)h(\tau)+g(\sigma)f(\tau)h(\sigma)\geq 0.
\end{equation}
\noindent Now multiplying equation (4.35) by $\mathrm{\psi}(x,\tau)$, which is positive because $\tau \in [1,x]$, $x>1$, then integrate the resulting identity with respect to $\tau$ from $1$ to $x$, we get
\begin{equation}
\begin{split}
&\frac{f(\sigma)}{\beta^{\alpha}\Gamma(\alpha)} \int_{1}^{x}e^{[\frac{\beta-1}{\beta}(ln x-ln \tau)]}(ln x-ln \tau)^{\alpha-1}[g(\tau) h(\tau)]\frac{d\tau}{\tau} \\
&-\frac{h(\sigma)}{\beta^{\alpha}\Gamma(\alpha)} \int_{1}^{x}e^{[\frac{\beta-1}{\beta}(ln x-ln \tau)]}(ln x-ln \tau)^{\alpha-1}[f(\tau) g(\tau)]\frac{d\tau}{\tau}\\
&+\frac{f(\sigma)g(\sigma)}{\beta^{\alpha}\Gamma(\alpha)} \int_{1}^{x}e^{[\frac{\beta-1}{\beta}(ln x-ln \tau)]}(ln x-ln \tau)^{\alpha-1}[h(\tau)]\frac{d\tau}{\tau}\\
&\frac{g(\sigma)h(\sigma)}{\beta^{\alpha}\Gamma(\alpha)} \int_{1}^{x}e^{[\frac{\beta-1}{\beta}(ln x-ln \tau)]}(ln x-ln \tau)^{\alpha-1}[f(\tau)]\frac{d\tau}{\tau} \geq 0,
\end{split}
\end{equation}
\noindent we get
\begin{equation}
\begin{split}
&f(\sigma)\mathcal{H}^{\alpha,\beta}_{1,x}[(gh)(x)]+g(\sigma)h(\sigma)\mathcal{H}^{\alpha,\beta}_{1,x}[f(x)]\\-
&g(\sigma)f(\rho)\mathcal{H}^{\alpha,\beta}_{1,x}[h(x)]-h(\sigma)\mathcal{H}^{\alpha,\beta}_{1,x}[(gf)(x)]\geq 0.
\end{split}
\end{equation}
\noindent Again multiplying (4.37)  by $\mathrm{\psi}(x,\sigma)$,in view of equation (3.3)and which is positive because $\sigma \in [1,x]$, $x>1$, then integrate the resulting identity with respect to $\sigma$ from $1$ to $x$, we get
\begin{equation}
\begin{split}
\mathcal{H}^{\alpha,\beta}_{1,x}[f(x)]\mathcal{H}^{\alpha,\beta}_{1,x}[(gh)(x)]&-\mathcal{H}^{\alpha,\beta}_{1,x}[h(x)]\mathcal{H}^{\alpha,\beta}_{1,x}[(gf)(x)]\\
&-\mathcal{H}^{\alpha,\beta}_{1,x}[(gf)(x)]\mathcal{H}^{\alpha,\beta}_{1,x}[h(x)]\\
 &+\mathcal{H}^{\alpha,\beta}_{1,x}[(gh)(x)]\mathcal{H}^{\alpha,\beta}_{1,x}[f(x)]\geq0,
\end{split}
\end{equation}
\noindent which implies that,
\begin{equation}
\mathcal{H}^{\alpha,\beta}_{1,x}[f(x)]\mathcal{H}^{\alpha,\beta}_{1,x}[(g h)(x)]\geq \mathcal{H}^{\alpha,\beta}_{1,x}[h(x)]\mathcal{H}^{\alpha,\beta}_{1,x}[(g f)(x)],
\end{equation}
\noindent we get
\begin{equation}
\frac{\mathcal{H}^{\alpha,\beta}_{1,x}[f(x)]}{\mathcal{H}^{\alpha,\beta}_{1,x}[h(x)]}\geq \frac{\mathcal{H}^{\alpha,\beta}_{1,x}[(gf)(x)]}{\mathcal{H}^{\alpha,\beta}_{1,x}[(g h)(x)]}.
\end{equation}
This completes the proof.
\begin{theorem} Suppose that $f$, $g$ and $h$ be positive and continuous functions on $[1,\infty)$, such that
 \begin{equation}
(g(\tau)-g(\sigma)) \left(\frac{f(\sigma)}{h(\sigma)}-\frac{f(\tau)}{h(\tau)}\right)\geq 0, \tau, \sigma \in(1,x)\ \  x>0,
\end{equation}
 then for all $\alpha>0,$ $\beta \in (0,1]$ $\alpha \in \mathbb{C}$ and $\mathfrak{R(\alpha)}> 0$,

\begin{equation}
\frac{\mathcal{H}^{\alpha,\beta}_{1,x}[f(x)]\mathcal{H}^{\phi,\varphi}_{1,x}[(gh)(x)]+\mathcal{H}^{\phi,\varphi}_{1,x}[f(x)]\mathcal{H}^{\alpha,\beta}_{1,x}[(gh)(x)]}{\mathcal{H}^{\alpha,\beta}_{1,x}[h(x)]\mathcal{H}^{\phi,\varphi}_{1,x}[(gf)(x)]+\mathcal{H}^{\phi,\varphi}_{1,x}[h(x)]\mathcal{H}^{\alpha,\beta}_{1,x}[(gf)(x)]}\geq 1,
\end{equation}
hold.
\end{theorem}

\textbf{Proof}:Multiplying equation (4.37) by $\frac{1}{\varphi^{\phi}\Gamma(\phi)\sigma}e^{[\frac{\varphi-1}{\varphi}(ln x-ln \sigma)]}(ln x-ln \sigma)^{\phi-1}$, which remain positive because$\sigma\in (1, x)$ , $(x>1), \phi, \varphi>0$. Then integrating resulting identity with respect to $\sigma$ from $1$ to $x$, we have,
\begin{equation}
\begin{split}
\mathcal{H}^{\phi,\varphi}_{1,x}[f(x)]\mathcal{H}^{\alpha,\beta}_{1,x}[(gh)(x)]&-\mathcal{H}^{\phi,\varphi}_{1,x}[h(x)]\mathcal{H}^{\alpha,\beta}_{1,x}[(gf)(x)]\\
&-\mathcal{H}^{\phi,\varphi}_{1,x}[(gf)(x)]\mathcal{H}^{\alpha,\beta}_{1,x}[h(x)]\\
 &+\mathcal{H}^{\phi,\varphi}_{1,x}[(gh)(x)]\mathcal{H}^{\alpha,\beta}_{1,x}[f(x)]\geq0,
\end{split}
\end{equation}
we get,
\begin{equation}
\begin{split}
&\mathcal{H}^{\phi,\varphi}_{1,x}[f(x)]\mathcal{H}^{\alpha,\beta}_{1,x}[(gh)(x)]+\mathcal{H}^{\phi,\varphi}_{1,x}[(gh)(x)]\mathcal{H}^{\alpha,\beta}_{1,x}[f(x)]\\
 &\geq \mathcal{H}^{\phi,\varphi}_{1,x}[h(x)]\mathcal{H}^{\alpha,\beta}_{1,x}[(gf)(x)]+\mathcal{H}^{\phi,\varphi}_{1,x}[(gf)(x)]\mathcal{H}^{\alpha,\beta}_{1,x}[h(x)],
 \end{split}
\end{equation}
this gives the required inequality (4.42).
\begin{theorem} Suppose that $f$ and $h$ are two positive continuous function such that $f\leq h$ on $[0, \infty)$. If $\frac{f}{h}$ is decreasing and $f$ is increasing on $[1, \infty)$, then for any $p\geq0$, For all $x >1$, $\alpha, \beta > 0$ we have
 \begin{equation}
\frac{\mathcal{H}^{\alpha,\beta}_{1,x}[f(x)]}{\mathcal{H}^{\alpha,\beta}_{1,x}[h(x)]}\geq \frac{\mathcal{H}^{\alpha,\beta}_{1,x}[f^{p}(x)]}{\mathcal{H}^{\alpha,\beta}_{1,x}[h^{p}(x)]}.
\end{equation}
\end{theorem}
\textbf{Proof}: We take $g=f^{p-1}$ in theorem 4.5.
\begin{equation}
\frac{\mathcal{H}^{\alpha,\beta}_{1,x}[f(x)]}{\mathcal{H}^{\alpha,\beta}_{1,x}[h(x)]}\geq \frac{\mathcal{H}^{\alpha,\beta}_{1,x}[(ff^{p-1})(x)]}{\mathcal{H}^{\alpha,\beta}_{1,x}[(hf^{p-1})(x)]}.
\end{equation}
Since  $f\leq h$ on  $[1, \infty)$, then we can write,
\begin{equation}
hf^{p-1}\leq h^{p}.
\end{equation}
Multiplying equation (4.47) by$\mathrm{\psi}(x,\tau)$, which is positive because $\tau \in [1,x]$, $x>1$, then integrate the resulting identity with respect to $\tau$ from $1$ to $x$, we get
\begin{equation}
\begin{split}
& \frac{1}{\beta^{\alpha}\Gamma(\alpha)} \int_{1}^{x}e^{[\frac{\beta-1}{\beta}(ln x-ln \tau)]}(ln x-ln \tau)^{\alpha-1}[f^{p-1}h(\tau)]\frac{d\tau}{\tau} \\
&\leq\frac{1}{\beta^{\alpha}\Gamma(\alpha)} \int_{1}^{x}e^{[\frac{\beta-1}{\beta}(ln x-ln \tau)]}(ln x-ln \tau)^{\alpha-1}[h^{p}(\tau)]\frac{d\tau}{\tau},
\end{split}
\end{equation}
implies that
\begin{equation}
\mathcal{H}^{\alpha,\beta}_{1,x}[hf^{p-1}(x)]\leq \mathcal{H}^{\alpha,\beta}_{1,x}[h^{p}(x)],
\end{equation}
and so we have,
\begin{equation}
 \frac{\mathcal{H}^{\alpha,\beta}_{1,x}[(f f^{p-1})(x)]}{\mathcal{H}^{\alpha,\beta}_{1,x}[(h f^{p-1})(x)]}\geq \frac{\mathcal{H}^{\alpha,\beta}_{1,x}[f^{p}(x)]}{\mathcal{H}^{\alpha,\beta}_{1,x}[h^{p}(x)]},
\end{equation}
then from equation (4.46) and (4.50), we obtain (4.45).
\begin{theorem} Suppose that $f$ and $h$ are two positive continuous function such that $f\leq h$ on $[1, \infty)$. If $\frac{f}{h}$ is decreasing and $f$ is increasing on $[1, \infty)[$, then for any $p\geq1$, for all $x >0$, $\sigma, \beta, \phi, \varphi>0$
 \begin{equation}
\frac{\mathcal{H}^{\alpha,\beta}_{1,x}[f(x)]\mathcal{H}^{\phi,\varphi}_{1,x}h^{p}[(x)]+\mathcal{H}^{\phi,\varphi}_{1,x}[f(x)]\mathcal{H}^{\alpha,\beta}_{1,x}[h^{p}(x)]}{\mathcal{H}^{\alpha,\beta}_{1,x}[h(x)]\mathcal{H}^{\phi,\varphi}_{1,x}[f^{p}(x)]+\mathcal{H}^{\phi,\varphi}_{1,x}[h(x)]\mathcal{H}^{\alpha,\beta}_{1,x}[f^{p}(x)]}\geq 1.
\end{equation}
\end{theorem}
\textbf{Proof}: We take $g=f^{p-1}$ in theorem 4.6, then we obtain
 \begin{equation}
\frac{\mathcal{H}^{\alpha,\beta}_{1,x}[f(x)]\mathcal{H}^{\phi,\varphi}_{1,x}[hf^{p-1}(x)]+\mathcal{H}^{\phi,\varphi}_{1,x}[f(x)]\mathcal{H}^{\alpha,\beta}_{1,x}[hf^{p-1}(x)]}{\mathcal{H}^{\alpha,\beta}_{1,x}[h(x)]\mathcal{H}^{\phi,\varphi}_{1,x}[f^{p}(x)]+\mathcal{H}^{\phi,\varphi}_{1,x}[h(x)]\mathcal{H}^{\alpha,\beta}_{1,x}[f^{p}(x)]}\geq 1,
\end{equation}
then by hypothesis, $f\leq h$ on $[1,\infty)$, which implies that
  \begin{equation}
  hf^{p-1}\leq h^{p}.
\end{equation}
Now, multiplying both side of (4.53) by  $\frac{1}{\varphi^{\phi}\Gamma(\phi)\sigma}e^{[\frac{\varphi-1}{\varphi}(ln x-ln \sigma)]}(ln x-ln \sigma)^{\phi-1}$, which remain positive because$\sigma\in (1, x)$ , $(x>1), \phi, \varphi>0$. Then integrating resulting identity with respect to $\sigma$ from $1$ to $x$, we have,
 \begin{equation}
\mathcal{H}^{\phi,\varphi}_{1,x}[h f^{p-1}(x)] \leq \mathcal{H}^{\phi,\varphi}_{1,x}[h^{p}(x)],
\end{equation}
multiplying on both side of (4.50) by $\mathcal{H}^{\alpha,\beta}_{1,x}[f(x)]$, we obtain
\begin{equation}
\mathcal{H}^{\alpha,\beta}_{1,x}[f(x)]\mathcal{H}^{\phi,\varphi}_{1,x}[h f^{p-1}(x)] \leq \mathcal{H}^{\alpha,\beta}_{1,x}[f(x)]\mathcal{H}^{\phi,\varphi}_{1,x}[h^{p}(x)],
\end{equation}
hence by (4.49) and (4.55), we obtain
\begin{equation}
\begin{split}
&\mathcal{H}^{\alpha,\beta}_{1,x}[f(x)]\mathcal{H}^{\phi,\varphi}_{1,x}[h f^{p-1}(x)]+\mathcal{H}^{\phi,\varphi}_{1,x}[f(x)] \mathcal{H}^{\alpha,\beta}_{1,x}[h f^{p-1}(x)]\\
&\leq \mathcal{H}^{\alpha,\beta}_{1,x}[f(x)]\mathcal{H}^{\phi,\varphi}_{1,x}[h^{p}(x)]+\mathcal{H}^{\phi,\varphi}_{1,x}[f(x)] \mathcal{H}^{\alpha,\beta}_{1,x}[h^{p}(x)].
\end{split}
\end{equation}
By (4.52) and (4.56), we complete the proof of this theorem.
		
\end{document}